\documentclass{amsart}

\newtheorem{theorem}{Theorem}[section]
\newtheorem{lemma}[theorem]{Lemma}
\newtheorem{corollary}[theorem]{Corollary}

\theoremstyle{definition}

\theoremstyle{remark}

\numberwithin{equation}{section}

\newcommand{\abs}[1]{\lvert#1\rvert}

\DeclareSymbolFont{AMSb}{U}{msb}{m}{n}
\DeclareMathSymbol{\Z}{\mathalpha}{AMSb}{"5A}

\DeclareMathSymbol{\nmid}{\mathrel}{AMSb}{"2D}

\begin{document}
\newcommand{\beqs}{\begin{equation*}}
\newcommand{\eeqs}{\end{equation*}}
\newcommand{\beq}{\begin{equation}}
\newcommand{\eeq}{\end{equation}}
\newcommand\nutwid{\overset {\text{\lower 3pt\hbox{$\sim$}}}\nu}
\newcommand\Mtwid{\overset {\text{\lower 3pt\hbox{$\sim$}}}M}
\newcommand\ptwid{\overset {\text{\lower 3pt\hbox{$\sim$}}}p}
\newcommand\pitwid{\overset {\text{\lower 3pt\hbox{$\sim$}}}\pi}
\newcommand\bijone{\overset {1}\longrightarrow}                   
\newcommand\bijtwo{\overset {2}\longrightarrow}                   
\newcommand\pihat{\widehat{\pi}}
\newcommand\mymod[1]{(\mbox{mod}\ {#1})}
\newcommand\myto{\to}
\newcommand\srank{\mathrm{srank}}
\newcommand\pitc{\pi_{\mbox{$t$-core}}}
\newcommand\pitcb[1]{\pi_{\mbox{${#1}$-core}}}  
\newcommand\mpitcb[1]{\pi_{\mathrm{{#1}-core}}}  
\newcommand\lamseq{(\lambda_1, \lambda_2, \dots,\lambda_\nu)}
\newcommand\mylabel[1]{\label{#1}}
\newcommand\eqn[1]{(\ref{eq:#1})}
\newcommand\stc{{St-crank}}
\newcommand\mstc{\mbox{St-crank}}
\newcommand\mmstc{\mbox{\scriptsize\rm St-crank}}
\newcommand\tqr{{$2$-quotient-rank}}
\newcommand\mtqr{\mbox{$2$-quotient-rank}}
\newcommand\mmtqr{\mathrm{2-quotient-rank}}
\newcommand\fcc{{$5$-core-crank}}        
\newcommand\bgr{BG-rank} 
\newcommand\mbgr{\mbox{BG-rank}} 
\newcommand\nvec{(n_0, n_1, \dots, n_{t-1})}
\newcommand\parity{\mbox{par}}
\newcommand\pbar{\overline{p}}
\newcommand\Pbar{\overline{P}}
\newcommand\abar{\overline{a}}
\newcommand\RHS{\mbox{RHS}}
\newcommand\LHS{\mbox{LHS}}
\newcommand\pli{\pi_i^L}
\newcommand\poo{\pi_1^1}
\newcommand\poz{\pi_0^1}
\newcommand\tpoz{\pitwid_0^1}
\newcommand\plo{\pi_1^L}
\newcommand\plz{\pi_0^L}
\newcommand\tplz{\pitwid_0^L}
\newcommand\tplm{\pitwid_m^L}
\newcommand\Atwid{\overset{\text{\lower 3pt\hbox{$\sim$}}}A}

\title[Dissecting the Stanley partition function]{Dissecting the Stanley partition function}

\author{Alexander Berkovich}
\address{Department of Mathematics, University of Florida, Gainesville,
Florida 32611-8105}
\email{alexb@math.ufl.edu}          

\author{Frank G. Garvan}
\address{Department of Mathematics, University of Florida, Gainesville,
Florida 32611-8105}
\email{frank@math.ufl.edu}          

\subjclass[2000]{Primary 11P81, 11P82, 11P83; Secondary 05A17, 05A19}

\date{\today}  


\keywords{generating functions, Stanley's partitions, 
even/odd dissection, upper bounds, asymptotic formulas,
partition inequalities}

\begin{abstract}
Let $p(n)$ denote the number of unrestricted partitions of $n$. For
$i=0$, $2$, let $p_i(n)$ denote the number of partitions $\pi$ of $n$
such that ${\mathcal O}(\pi) - {\mathcal O}(\pi') \equiv i\pmod{4}$.
Here ${\mathcal O}(\pi)$ denotes the number of odd parts of the partition $\pi$
and $\pi'$ is the conjugate of $\pi$. 
R. Stanley \cite{St1}, \cite{St2} derived an infinite product 
representation for the generating function of $p_0(n)-p_2(n)$. 
Recently, H. Swisher \cite{Sw} employed the circle method to show that 
\beqs
\lim_{n\to\infty} \frac{p_0(n)}{p(n)} = \frac{1}{2},
\tag{i}
\eeqs
and that for sufficiently large $n$
\beqs
\begin{array}{ll}
2 p_0(n) > p(n), &\qquad\mbox{if $n\equiv0,1\pmod{4}$},\\
2 p_0(n) < p(n), &\qquad\mbox{otherwise}.
\end{array}
\tag{ii}
\eeqs

In this paper we study the even/odd dissection of the Stanley product,
and show how to use it to prove (i) and (ii) with no restriction on $n$.
Moreover, we establish the following new result
\beqs
\abs{p_0(2n) - p_2(2n)} > \abs{p_0(2n+1) - p_2(2n+1)}, \qquad n>0.
\eeqs
Two proofs of this surprising inequality are given. The first one uses the 
G\"ollnitz-Gordon partition theorem. The second one is an immediate corollary
of a new partition inequality, which we prove in a combinatorial manner.
Our methods are elementary. 
We use only Jacobi's triple product identity and some naive upper bound
estimates.
\end{abstract}

\maketitle

\section{Introduction} \label{sec:intro}
Let $\pi$ denote a partition of some integer and $\pi'$ its conjugate.
Let ${\mathcal O}(\pi)$ denote the number of odd parts of $\pi$ and
$p_i(n)$ denote the number of partitions of $n$ for which
${\mathcal O}(\pi) - {\mathcal O}(\pi') \equiv i\pmod{4}$. It is easy to see
that
\beq
n\equiv{\mathcal O}(\pi)\equiv{\mathcal O}(\pi') \pmod{2}
\mylabel{eq:1}
\eeq
for any partition $\pi$ of $n$, so that
\beq
p(n) = p_0(n) + p_2(n),
\mylabel{eq:2}
\eeq
where $p(n)$ is the number of unrestricted partitions of $n$.
Obviously,
\beq
\sum_{n\ge0}(p_0(n) + p_2(n))q^n = \sum_{n\ge0} p(n) q^n 
= \frac{1}{(q;q)_\infty}.
\mylabel{eq:3}
\eeq
We use the standard notation
\begin{align}
(a;q)_\infty &= \lim_{L\to\infty} (a;q)_L, \mylabel{eq:4} \\ 
(a;q)_L = (a)_L &=
\begin{cases}
1, &\mbox{if $L=0$},\\
\prod_{j=0}^{L-1}(1-aq^j), &\mbox{if $L>0$,}
\end{cases}
\mylabel{eq:4b} 
\end{align}
\beq
(a_1,a_2,\dots,a_n;q)_\infty = (a_1;q)_\infty (a_2;q)_\infty \cdots
(a_n;q)_\infty,
\mylabel{eq:5}
\eeq
and 
\beq
(a_1,a_2,\dots,a_n;q)_L = (a_1;q)_L (a_2;q)_L \cdots
(a_n;q)_L.
\mylabel{eq:5b}
\eeq
Recently, R. Stanley \cite{St1}, \cite{St2} has shown that
\beq
\sum_{n\ge0}(p_0(n) - p_2(n))q^n =
\frac{ (-q;q^2)_\infty}{(q^4,-q^2,-q^2;q^4)_\infty}.
\mylabel{eq:6}
\eeq
G. Andrews \cite{An1} used \eqn{3} and \eqn{6} to show that
\beq
\sum_{n\ge0} p_0(n) q^n = 
\frac{E^2(q^2) E^5(q^{16})}
{E(q) E^5(q^4) E^2(q^{32})},
\mylabel{eq:7}
\eeq
where $E(q) := (q;q)_\infty$.
Moreover, he proved that
\beq
p_0(5n+4) \equiv 0 \pmod{5},
\mylabel{eq:8}
\eeq
which is a refinement of the famous Ramanujan congruence \cite{R}
\beq
p(5n+4) \equiv 0 \pmod{5}.
\mylabel{eq:9}
\eeq
Various combinatorial proofs of \eqn{6} and its generalizations
were given by A. Sills \cite{Si1}, C. Boulet \cite{B}
and A.J. Yee \cite{Y}. A combinatorial proof of \eqn{8}
was found by A. Berkovich and F. Garvan \cite{BG}. In a recent paper
\cite{Sw}, H. Swisher showed that \eqn{8} is just one of infinitely many
similar congruences satisfied by $p_0(n)$. In addition, she applied
the Hardy-Ramanujan `circle method' \cite{HR}, \cite{An2} to the
product in \eqn{7} to deduce two interesting corollaries:
\begin{corollary}[Swisher]
\label{cor1}
\beqs
\lim_{n\to\infty} \frac{p_0(n)}{p(n)}=\frac{1}{2}.
\eeqs
\end{corollary}
\begin{corollary}[Swisher]
\label{cor2} If $n$ is sufficiently large, then
\begin{enumerate}
\item[(a)] $p_0(n) > \tfrac{1}{2} p(n)$, if $n\equiv 0,1 \pmod{4}$,
\item[(b)] $p_0(n) < \tfrac{1}{2} p(n)$, if $n\equiv 2,3 \pmod{4}$.
\end{enumerate}
\end{corollary}
One object of this paper is to provide elementary proofs of Corollary 
\ref{cor1}, and  
Corollary \ref{cor2} with the restriction ``$n$ is sufficiently large'' 
removed. 
To this end we will prove a Dissection Theorem for
the Stanley infinite product \eqn{6}:
\begin{theorem}[Dissection Theorem]
\label{disthm}
\beq
\frac{(-q;q^2)_\infty}
{(q^4,-q^2,-q^2;q^4)_\infty} = F_0(-q^2) + q F_1(-q^2),
\mylabel{eq:10}
\eeq
where
\beq
F_i(q) = \frac{1}{E(q)}
\frac{1}
{(q^{1+2i},q^2,q^4,q^6,q^{7-2i};q^8)_\infty},
\mylabel{eq:11}
\eeq
for $i=0,1$.
\end{theorem}
The proof of this theorem, given in the next section, requires only
the Jacobi triple product identity:
\beq
\sum_{n=-\infty}^\infty q^{n^2} z^n
= (q^2,-zq,-q/z;q^2)_\infty.
\mylabel{eq:12}
\eeq
We will show that the Dissection Theorem immediately implies 
Corollary \ref{cor2} with no restriction on $n$. In Section 3 we will use
only elementary methods to prove the following upper bound:
\begin{lemma}[Upper Bound Lemma]
\label{uplem}
For $n\ge0$ and $i=0,1$
\beq
|p_0(2n+i) - p_2(2n+i)| < \exp\{\tfrac{\pi}{2} \sqrt{\tfrac{13n}{3}}\}.
\mylabel{eq:13}
\eeq
\end{lemma}
Hardy and Ramanujan \cite{HR} established in their classical paper
that
\beq
p(n) \sim \frac{A}{n} \exp\{ \pi \sqrt{\tfrac{2n}{3}}\},
\mylabel{eq:14}
\eeq
with $A=\frac{1}{4\sqrt{3}}$.
An elementary proof of \eqn{14} (with undetermined $A$) was
given later by Erd{\H o}s \cite{E}. Obviously,
\eqn{13} and \eqn{14} along with \eqn{2} imply Corollary \ref{cor1}.
In Section 4 we will sharpen the upper bound in \eqn{13} and 
prove the following new result
\beq
\abs{p_0(2n) - p_2(2n)} > \abs{p_0(2n+1) - p_2(2n+1)}, \qquad n\ge1.
\mylabel{eq:14b}
\eeq
Our first proof of \eqn{14b} makes use of a relation between $F_0(q)$,
$F_1(q)$ and the G\"ollnitz-Gordon products.  Also we show, using Meinardus's Theorem,
that
\beq
\lim_{n\to\infty} \frac{p_0(2n) - p_2(2n)}{p_0(2n+1)-p_2(2n+1)}
= 1 + \sqrt{2}.
\mylabel{eq:14c}
\eeq
In Section 5, we will
establish a new partition inequality from which \eqn{14b} follows as
an easy corollary. We conclude with some conjectures.

\section{Proof of the Dissection Theorem and the strong version of Corollary \ref{cor2}}
\label{sec:2}
We begin by observing that
\beq
E(q) = (q;q)_\infty = (q,q^2;q^2)_\infty,
\mylabel{eq:15}
\eeq
and
so
\beq
E(-q) = (-q,q^2;q^2)_\infty.
\mylabel{eq:16}
\eeq
This allows us to rewrite the right side of \eqn{10} as

\begin{align}
\RHS\eqn{10} &= 
\frac{1}{(-q^2,q^4;q^4)_\infty}
\left\{
\frac{1}{(-q^2,q^4,q^8,q^{12},-q^{14};q^{16})_\infty}  \right.
\mylabel{eq:16b}\\
&\qquad\qquad\qquad \left.
 + \frac{q}{(-q^6,q^4,q^8,q^{12},-q^{10};q^{16})_\infty} \right\}
\nonumber \\
&= \frac{E(q^{16})}{(-q^2,q^4,q^4;q^4)_\infty}
\left\{
\frac{1}{(-q^2,-q^{14};q^{16})_\infty} +
\frac{q}{(-q^6,-q^{10};q^{16})_\infty}
\right\},
\nonumber
\end{align}
where we have used 
\beq
(q^4;q^4)_\infty = (q^4,q^8,q^{12},q^{16};q^{16})_\infty.
\mylabel{eq:17}
\eeq
Next, we employ
\begin{gather}
(-q^2,-q^6,-q^{10},-q^{14};q^{16})_\infty = (-q^2;q^4)_\infty,
\mylabel{eq:18}\\
(-q,-q^3;q^4)_\infty = (-q;q^2)_\infty,
\mylabel{eq:19}
\end{gather}
and Jacobi's triple product identity \eqn{12} to obtain
\begin{align}
\RHS\eqn{10} &=
\frac{
(q^{16},-q^6,-q^{10};q^{16})_\infty 
+ q (q^{16},-q^2,-q^{14};q^{16})_\infty }
{ (-q^2,q^4;q^4)_\infty^2 } \mylabel{eq:20}\\
&= 
\frac{
\sum_{j=-\infty}^\infty q^{8j^2+2j} + q \sum_{j=-\infty}^\infty q^{8j^2-6j}
}
{(-q^2,q^4;q^4)_\infty^2 }
=
\frac{\sum_{j=-\infty}^\infty q^{2j^2+j}}
{(-q^2,q^4;q^4)_\infty^2 } \nonumber\\
&= 
\frac{ (q^4,-q,-q^3;q^4)_\infty}
{(-q^2,q^4;q^4)_\infty^2 } 
=
\frac{ (-q;q^2)_\infty}
{(q^4,-q^2,-q^2;q^4)_\infty } = \LHS\eqn{10},
\nonumber
\end{align}
as asserted.

Before we move on we would like to point out that K. Alladi \cite{Al} studied
even/odd splits of many classical series.
In particular, he treated the Euler pentagonal series, the Gauss triangular
series and the famous Rogers-Ramanujan series.

It follows from the Dissection Theorem that for $i=0,1$
\beq
\sum_{n\ge0}(p_0(2n+i) - p_2(2n+i))q^n
= \frac{1}{E(-q)}
\frac{1}
{(-q^{1+2i},q^2,q^4,q^6,-q^{7-2i};q^8)_\infty}.
\mylabel{eq:21}
\eeq
Replacing $q$ by $-q$ in \eqn{21} we find that
\beq
\sum_{n\ge0} (-1)^n(p_0(2n+i) - p_2(2n+i))q^n
= \frac{1}{E(q)}
\frac{1}
{(q^{1+2i},q^2,q^4,q^6,q^{7-2i};q^8)_\infty},
\mylabel{eq:22}
\eeq
for $i=0,1$. It is now obvious that for $n\ge0$ and $i=0,1$
\beq
(-1)^n(p_0(2n+i) - p_2(2n+i)) > 0.
\mylabel{eq:23}
\eeq
Recalling \eqn{2}, we see that  for $n\ge0$ and $i=0,1$
\beq
(-1)^n(2 p_0(2n+i) - p(2n+i)) > 0.
\mylabel{eq:24}
\eeq
In other words, we have the following corollary: For $n\ge0$,
\begin{enumerate}
\item[(a)]
$p_0(n) > \frac{p(n)}{2}$, if $n\equiv0,1\pmod{4}$,
\item[(b)]
$p_0(n) < \frac{p(n)}{2}$, if $n\equiv2,3\pmod{4}$.
\end{enumerate}
This corollary obviously implies Corollary \ref{cor2}.

\section{Proof of the Upper Bound Lemma and Corollary \ref{cor1}}
\label{sec:3}
Let $c_n$ denote $\abs{p_0(n) - p_2(n)}$.
Obviously,
\beq
c_n >0, \qquad n\ge0,
\mylabel{eq:25}
\eeq
and
\beq
\sum_{n\ge0} c_{2n+i} q^n = F_i(q), \qquad i=0,1.
\mylabel{eq:26}
\eeq
To obtain the upper bound \eqn{13} for the $c_n$ we will employ
the standard elementary argument \cite[pp.316--318]{Ap}, \cite{Od}. Assume $0<q<1$ so that
\beq
c_{2n+i} q^n < F_i(q), \quad \mbox{$i=0,1$ and $n\ge0$}.
\mylabel{eq:27}
\eeq
Clearly, \eqn{27} is a simple consequence of \eqn{25}, \eqn{26}, and
so we have for $n\ge0$, $i=0,1$               
\beq
\log(c_{2n+i}) < \log F_i(q) + ns,
\mylabel{eq:28}
\eeq
where $q=e^{-s}$ and $s>0$. To proceed further we make use of
\begin{align}
-\log(1-x) &= \sum_{m\ge1} \frac{x^m}{m},
\mylabel{eq:29}\\
\frac{1}{1-x} &= \sum_{n\ge0}x^n,
\mylabel{eq:30}
\end{align}
to find that for $i=0,1$
\beq
\log F_i(q) = \sum_{m\ge1} \frac{1}{m} \frac{1}{e^{sm}-1}
+ \sum_{\substack{m\ge 1\\ r\in S_i}} \frac{1}{m} \frac{e^{rms}}{e^{8sm}-1},
\mylabel{eq:31}
\eeq
where $S_i=\{1+2i,2,4,6,7-2i\}$. Next, we shall require the
following inequalities:
\beq
\frac{1}{e^x-1} < \frac{1}{x}, \quad\mbox{$x>0$},
\mylabel{eq:32}
\eeq
\beq
\frac{e^{rx}+e^{(8-r)x}}{e^{8x}-1} < \frac{1}{4x}, \quad
\mbox{$x>0$, $r=2,3,4$,}
\mylabel{eq:33}
\eeq
and
\beq
\frac{e^x+e^{7x}+e^{4x}+e^{2x}+e^{6x}}{e^{8x}-1} < \frac{5}{8x}, \quad
\mbox{$x>0$}.
\mylabel{eq:34}
\eeq
We will prove \eqn{32}--\eqn{34} later. In the mean time we observe
that these inequalities imply that
\beq
\log F_i(q) < \frac{1}{s} \sum_{m\ge1}\frac{1}{m^2} +
\frac{5}{8s} \sum_{m\ge1}\frac{1}{m^2} 
= \frac{13}{48} \frac{\pi^2}{s}.
\mylabel{eq:35}
\eeq
Combining \eqn{28} and \eqn{35} we have for $i=0,1$ and $s>0$ that
\beq
\log c_{2n+i} < \frac{13}{48} \frac{\pi^2}{s} + ns.
\mylabel{eq:36}
\eeq
To minimize the right side of \eqn{36} we choose $s=\pi\sqrt{13/(48n)}$ to
find
\beq
\log c_{2n+i} < \frac{\pi}{2} \sqrt{\frac{13n}{3}}
\mylabel{eq:37}
\eeq
for $i=0,1$ and $n\ge0$. The above inequality \eqn{37} is essentially
\eqn{13}, as desired.

Obviously, \eqn{13} and \eqn{14} imply that
\beq
\lim_{n\to\infty} \frac{c_n}{p(n)} = 0.
\mylabel{eq:38}
\eeq
Next, using \eqn{2} and
\beq
p_0(2n+i)-p_2(2n+i) = c_{2n+i} (-1)^n,
\mylabel{eq:39}
\eeq
we get
\beq
\frac{p_0(2n+i)}{p(2n+i)} = \frac{1}{2} + (-1)^n \frac{c_{2n+i}}{p(2n+i)}.
\mylabel{eq:40}
\eeq
Corollary \ref{cor1} follows easily from \eqn{38} and \eqn{40}.
To complete the proof all we need is to verify \eqn{32}--\eqn{34}. 
To this end we recall that
\beq
e^x = \sum_{n\ge0} \frac{x^n}{n!},
\mylabel{eq:41}
\eeq
and that
\beq
e^x-1>0,\qquad\mbox{if $x>0$}.
\mylabel{eq:42}
\eeq
This allows us to rewrite \eqn{32} as the obvious relation
\beq
1 + x < e^x, \qquad\mbox{if $x>0$},
\mylabel{eq:43}
\eeq
and \eqn{33} as
\beq
4x\left( e^{(4-r)x} + e^{-(4-r)x}\right) < e^{4x} - e^{-4x},
\qquad\mbox{if $x>0$}
\mylabel{eq:44}
\eeq
with $r=2,3,4$. 
For $r=2,4$ \eqn{44} can be reduced to the obvious relation
\beq
rx < \frac{e^{rx}-e^{-rx}}{2}, \qquad\mbox{if $x>0$}.
\mylabel{eq:45}
\eeq
For $r=3$ it is equivalent to
\beq
2x < \frac{e^{3x}-e^{-3x}}{2} - \frac{e^{x}-e^{-x}}{2}, 
\qquad\mbox{if $x>0$},
\mylabel{eq:46}
\eeq
which follows from
\beq
0 < \sum_{n>0} \frac{3^{2n+1}-1}{(2n+1)!}x^{2n+1}, \qquad\mbox{if $x>0$}.
\mylabel{eq:47}
\eeq
To prove \eqn{34} we rewrite it as
\beq
8x \left( \frac{e^{3x}+e^{-3x}}{2} + \frac{e^{2x}+e^{-2x}}{2} + \frac{1}{2}
\right) < 5 \frac{e^{4x}-e^{-4x}}{2}.
\mylabel{eq:48}
\eeq
Using \eqn{41} we can reduce it to
\beq
8x \sum_{n>0} \frac{9^n+4^n}{(2n)!} x^{2n}
< 5 \sum_{n>0} \frac{(4x)^{2n+1}}{(2n+1)!}, 
\qquad\mbox{if $x>0$},
\mylabel{eq:49}
\eeq
which follows from
\beq
8 (2n+1) (9^n + 4^n) < 20 (16)^n, \qquad(n\ge1).
\mylabel{eq:50}
\eeq
Finally, \eqn{50} can be easily proven by a straightforward induction argument.

\section{Further Observations}
\label{sec:4}
It is possible to sharpen the upper bound in \eqn{13} with a little more effort.
To this end we note that for $i=0,1$
\begin{align}
& \sum_{n\ge0} (c_{2(n+1)+i} - c_{2n+i}) q^{n+1} + c_i = (1 - q) F_i(q)
\mylabel{eq:51}\\
& \quad = \frac{1}{(q^2;q)_\infty} 
\frac{1}{(q^{1+2i},q^2,q^4,q^6,q^{7-2i};q^8)_\infty}.
\nonumber
\end{align}
This means that for $n\ge0$
\beq
c_{2(n+1)+i} - c_{2n+i} \ge 0.
\mylabel{eq:52}
\eeq
Again we assume $0 < q < 1$.
So, instead of \eqn{27}, we have the inequality
\beq
\frac{c_{2n+i} q^n}{1-q} \le \sum_{k\ge n} c_{2k+i} q^k \le F_i(q).
\mylabel{eq:53}
\eeq
Letting $q=e^{-s}$ and taking logarithms we obtain
\beq
\log c_{2n+i} \le F_i(q) + ns + \log(1-e^{-s}),
\mylabel{eq:54}
\eeq
for $i=0,1$, $n\ge0$ and $s>0$. Moreover, for $s>0$ we have
\beq
\log(1-e^{-s}) < \log s.
\mylabel{eq:55}
\eeq
It follows from \eqn{35}, \eqn{54}, \eqn{55} that
\beq
\log c_{2n+i} < \frac{13\pi^2}{48s} + ns + \log s,
\mylabel{eq:56}
\eeq
for $i=0,1$, $n\ge0$ and $s>0$.
Evaluating the right side of \eqn{56} at
\beqs
s = \frac{\pi}{2} \sqrt{\frac{13}{12}} \frac{1}{\sqrt{n}},
\eeqs
we get, instead of \eqn{13}, the sharper upper bound
\beq
c_{2n+i} < \frac{\pi}{2} \sqrt{\frac{13}{12}} \frac{1}{\sqrt{n}}
 \exp\{\tfrac{\pi}{2} \sqrt{\tfrac{13n}{3}}\},
\mylabel{eq:57}
\eeq
for $n\ge0$ and $i=0,1$.
It would be difficult to improve on \eqn{57} using only elementary
methods. However,
applying Meinardus's Theorem \cite{M} (see also Th. 6.1 in \cite{An2})
to the products $F_0(q)$, $F_1(q)$, we obtain
\beq
c_{2n+i} \sim  \frac{\sqrt{\frac{13}{6}}}{32\sin((2i+1)\pi/8)}
\frac{1}{n}
 \exp\{\tfrac{\pi}{2} \sqrt{\tfrac{13n}{3}}\}\qquad (i=0,1),
\mylabel{eq:58}
\eeq
as $n\to\infty$.
It follows that
\beq
\lim_{n\to\infty} \frac{p_0(2n) - p_2(2n)}{p_0(2n+1)-p_2(2n+1)}
= \cot \tfrac{\pi}{8} = 1 + \sqrt{2}.
\mylabel{eq:58b}
\eeq
Also, it is clear that for sufficiently large $n$
\beq
c_{2n} > c_{2n+1}.
\mylabel{eq:59}
\eeq
Remarkably, \eqn{59} holds for all $n\ge1$.
In order to prove this, we note that
\beq
\sum_{n\ge0} \left(c_{2n}-c_{2n+1}\right) q^n = F_0(q) - F_1(q)
= \frac{1}{E(q)(q^2;q^4)_\infty}
\left\{G_0(q) - G_1(q)\right\},
\mylabel{eq:60}
\eeq
where
\beq
G_i(q) := \frac{1}{(q^{1+2i},q^4,q^{7-2i};q^8)_\infty}.
\mylabel{eq:61}
\eeq
According to the G\"ollnitz-Gordon partition theorem
(see Th. 7.11 with $k=2$ in \cite{An2}) $G_i(q)$, with $i=0,1$,
is the generating function for partitions into parts differing by at least
$2$ and having no consecutive even parts. In addition, at most $1-i$
parts are $\le2$. It is now clear that the coefficients in the expansion
\beq
G_0(q) - G_1(q) = \sum_{k\ge0} b_k q^k
\mylabel{eq:62}
\eeq
are all nonnegative. It is easy to check that
\beq
b_0=b_3=0, \quad b_1=b_2=1.
\mylabel{eq:63}
\eeq
Moreover, for $k\ge4$
\beq
b_k > 0.
\mylabel{eq:64}
\eeq
This is because for $k\ge4$ there is at least one partition, 
namely $1+(k-1)$, which is generated by $G_0(q)$ but not by $G_1(q)$.
Next, it is obvious that
\beq
d_k>0,
\mylabel{eq:65}
\eeq
for $k\ge0$. Here the $d_k$ are given by
\beq
\frac{1}{E(q)(q^2;q^4)_\infty} = \sum_{k\ge0} d_k q^k.
\mylabel{eq:66}
\eeq
It follows from \eqn{60}, \eqn{63}, \eqn{64}, \eqn{65} that
\beq
c_{2n} > c_{2n+1},\quad\mbox{if $n\ge1$.}
\mylabel{eq:67}
\eeq
In other words,
\beq
\abs{p_0(2n) - p_2(2n)} > \abs{p_0(2n+1) - p_2(2n+1)},\quad\mbox{if $n\ge1$,}
\mylabel{eq:68}
\eeq
as asserted earlier.

We note that \eqn{67} can be proven directly without any appeal to the
G\"ollnitz-Gordon partition theorem. In fact, all we need is to show that
\beq
e_j\ge0,\quad\mbox{if $j\ge2$,}
\mylabel{eq:69}
\eeq
where the $e_j$ are given in the expansion
\beq
\frac{1}{(q,q^7;q^8)_\infty}
-
\frac{1}{(q^3,q^5;q^8)_\infty}
= q + \sum_{j\ge2} e_j q^j.
\mylabel{eq:70}
\eeq
In the next section we will establish a new partition inequality, which
immediately implies \eqn{69}.

\section{A Partition Inequality}
\label{sec:5}
\begin{theorem}
\label{thm2}   
Let $A_{L,i}(n)$ denote the number of partitions of $n$ into parts
$\equiv\pm (1+2i)\pmod{8}$, such that
the largest part $\le 8L-2i-1$. Then
\beq
A_{L,0}(n) \ge A_{L,1}(n),
\mylabel{eq:71}
\eeq
where inequality is strict for $L\ge1$, $n\ne 0$, $3$, $5$, $6$.
\end{theorem}
Obviously, the generating function for $A_{L,i}(n)$ is given by
\beq
\sum_{n\ge0} A_{L,i}(n) q^n = \frac{1}{(q^{1+2i},q^{7-2i};q^{8})_L},
\mylabel{eq:72}
\eeq
and so \eqn{69} follows from \eqn{71} in the limit as $L\to\infty$.

To proceed further we shall require the following

\vskip 10pt
\noindent
\underbar{Notation.} Let $\abs{\pi}$ denote the norm (sum of parts) of 
a partition $\pi$. Let $\nu(\pi,i)$ and $\mu(\pi,i)$ denote
the number of parts of $\pi$ congruent to $i\pmod{8}$ and equal to $i$,
respectively. Let $\pli$ denote some partition generated by \eqn{72}.

We are now ready to prove \eqn{71} for $L=1$. 
We consider $\poo$. To define a corresponding partition $\poz$
we consider three cases. 

\vskip 10pt
\noindent
\underbar{Case 1:} $\mu(\poo,3)\ge\mu(\poo,5)$.
Obviously, $\poo$ consists of $\mu(\poo,5)$ pairs of the form $(3+5)$
and of $\mu(\poo,3)-\mu(\poo,5)$ unpaired $3$'s.
Let us rewrite each pair as $(1+7)$ and each unpaired $3$ as $(1+1+1)$. In this
way we obtain a partition $\poz$ such that $\abs{\poo}=\abs{\poz}$ and, 
in addition,
$\mu(\poz,1)\ge\mu(\poz,7)$, and $3\mid (\mu(\poz,1)-\mu(\poz,7))$.

\vskip 10pt
\noindent
\underbar{Case 2:} $\mu(\poo,5)>\mu(\poo,3)$, and 
$3\nmid (\mu(\poo,5)-\mu(\poo,3))$.
This time we have $\mu(\poo,3)$ pairs of the form $(3+5)$
and $\mu(\poo,5)-\mu(\poo,3)$ unpaired $5$'s.
As before, we rewrite each pair as $(1+7)$.
We replace each unpaired $5$ by $(1+1+1+1+1)$. In this way
we create a partition $\poz$ such that $\abs{\poo}=\abs{\poz}$ and,
in addition,
$\mu(\poz,1)>\mu(\poz,7)$, $5\mid (\mu(\poz,1)-\mu(\poz,7))$, and
$3\nmid (\mu(\poz,1)-\mu(\poz,7))$.

\vskip 10pt
\noindent
\underbar{Case 3:} $\mu(\poo,5)>\mu(\poo,3)$, and 
$3\mid (\mu(\poo,5)-\mu(\poo,3))$.
Again, we rewrite each pair $(3+5)$ as $(1+7)$, and each but the last two
unpaired $5$'s as $(1+1+1+1+1)$. The last two $5$'s we replace
by $7+1+1+1$. This way we obtain a partition $\poz$ such that 
$\abs{\poo}=\abs{\poz}$ and,
in addition,
$\mu(\poz,1)>\mu(\poz,7)>0$, and $\mu(\poz,1)-\mu(\poz,7)\equiv7\pmod{15}$.

Clearly, cases 1, 2, and 3 above describe a map from the partitions
generated by $\frac{1}{(1-q^3)(1-q^5)}$
to the partitions generated by $\frac{1}{(1-q)(1-q^7)}$.
It is important to observe that this map is $1$--$1$ but not onto.
This means that $A_{1,0}(n)\ge A_{1,1}(n)$, for all $n\ge0$.
We show that this inequality is strict for $n\neq 0,3,5,6$, by constructing
a partition counted by $A_{1,0}(n)$ but which is not the image
of some partition counted by $A_{1,1}(n)$. First, we observe
that in each case $\poz$ satisfies
\beq
\mu(\poz,1)\ge\mu(\poz,7).
\mylabel{eq:73}
\eeq
For $m>r$ the partition $(7^m,1^r)$ does not satisfy \eqn{73},
and so this is the desired partition for the case $n=7m+r$, $m>r$
and $r=0$, $1$, \dots, $6$. 
To complete the proof of our assertion we need to examine integers, which are not
of the form $7m+r$ for some $m>r$ and $0\le r<7$.
It is easy to check that all such integers form the set $\tilde S_1\bigcup \tilde S_2\bigcup \tilde S_3$.
Here,
\begin{align*}
\tilde S_1:=&\{1,2,4,8,11,13,16,17,19,26,32,34,41\},\\
\tilde S_2:=&\{9,10,12,24,25,27,40\},\\
\tilde S_3:=&\{18,20,33,48\}.
\end{align*}
If $n\in \tilde S_1$, then the desired partition is $(1^n)$.
If $n\in \tilde S_2$, then the desired partition is $(7^1,1^{(n-7)})$.
Finally, if $n\in \tilde S_3$, then the desired partitions are
$(7^2,1^4)$, $(7^2,1^6)$, $(7^4,1^5)$, $(7^4,1^{20})$.
To see that this is indeed the case we observe that all constructed partitions satisfy
$3\nmid(\mu(\pi^1_0,1)-\mu(\pi^1_0,7))$, 
$5\nmid(\mu(\pi^1_0,1)-\mu(\pi^1_0,7))$, 
$(\mu(\pi^1_0,1)-\mu(\pi^1_0,7))\not\equiv 7\pmod{15}$.

It remains to prove \eqn{71} for $L>1$. We start by removing the 
multiples of $8$ from each part of $\plo$. Next, we
assemble the extracted multiples of $8$ from the parts congruent
to $3,5\pmod{8}$ into two vectors $8\vec{v}_3$ and $8\vec{v}_5$.
The vectors $\vec{v}_3$, $\vec{v}_5$ have nonnegative
integer components, arranged in nondecreasing order. The dimensions of these
vectors are $\nu(\plo,3)$ and $\nu(\plo,5)$, respectively.
Having extracted the multiples of $8$ from the parts of $\plo$, we 
obtain a partition $\poo$. Next we convert $\poo$ into $\poz$ using the
map described above. Then we need to reattach the multiples of $8$ to
the parts of $\poz$. The procedure depends on the same three cases we 
considered earlier. 

\vskip 10pt
\noindent
\underbar{Case 1.}
We add the components of $8\vec{v}_3$,
$8\vec{v}_5$ to the parts of $\poz$ that are equal to $1$, $7$,
respectively. This way we create a partition $\plz$
that satisfies 
$\nu(\plz,1)\ge\nu(\plz,7)$, $3\mid (\nu(\plz,1)-\nu(\plz,7))$, and
$\mu(\plz,1) \ge \tfrac{2}{3}(\nu(\plz,1)-\nu(\plz,7))$.
To understand this inequality observe that
$\mu(\plz,1)\ge\nu(\plz,1)-\nu(\pi^L_1,3)=
3(\nu(\pi^L_1,3)-\nu(\pi^L_1,5))+\nu(\pi^L_1,5)-\nu(\pi^L_1,3)=
2(\nu(\pi^L_1,3)-\nu(\pi^L_1,5))=\tfrac{2}{3}(\nu(\pi^L_0,1)-\nu(\pi^L_0,7))$.
And so, $\mu(\pi^L_0,1)\ge\tfrac{2}{3}(\nu(\pi^L_0,1)-\nu(\pi^L_0,7))$, as claimed.

\vskip 10pt
\noindent
\underbar{Case 2.}
We add the components of $8\vec{v}_3$,
$8\vec{v}_5$ to the parts of $\poz$ that are equal to $7$, $1$,
respectively. This way we obtain a partition $\plz$
such that  
$\nu(\plz,1)\ge\nu(\plz,7)$, $5\mid (\nu(\plz,1)-\nu(\plz,7))$, 
$3\nmid (\nu(\plz,1)-\nu(\plz,7))$, and 
$\mu(\plz,1) \ge \tfrac{4}{5}(\nu(\plz,1)-\nu(\plz,7))$.

\vskip 10pt
\noindent
\underbar{Case 3.}
As in Case 2, we add the components of $8\vec{v}_3$,
$8\vec{v}_5$ to the parts of $\poz$ that are equal to $7$, $1$,
respectively, and  obtain a partition $\plz$
such that 
$\nu(\poz,1)>\nu(\poz,7)$, $\nu(\poz,1)-\nu(\poz,7)\equiv7\pmod{15}$, 
$\mu(\plz,7)>0$ and
$\mu(\plz,1) \ge \tfrac{1}{5}(4(\nu(\plz,1)-\nu(\plz,7))-3)$.

We illustrate our map with the following example.
Let $\pi^{11}_1=(3,19,43,45^2,53,85)$ be a partition of $293$.
Note that $\nu(\pi^{11}_1,3)=3$, $\mu(\pi^{11}_1,3)=1$, $\nu(\pi^{11}_1,5)=4$, $\mu(\pi^{11}_1,5)=0$.
This partition gives rise to $\vec{v}_3=(0,2,5)$, $\vec{v}_5=(5,5,6,10)$,    
$\pi^1_1=(3^3,5^4)$, $\pi^1_0=(1^8,7^3)$, $\pi^{11}_0=(1^4,7,23,41^2,47,49,81)$.
Note that $\nu(\pi^{11}_0,1)=8$, $\mu(\pi^{11}_0,1)=4$, $\nu(\pi^{11}_0,7)=3$,  $\mu(\pi^{11}_0,7)=1$.
And so,
$\nu(\pi^{11}_0,1)>\nu(\pi^{11}_0,7)$, $5\mid(\nu(\pi^{11}_0,1)-\nu(\pi^{11}_0,7))$,
$3\nmid(\nu(\pi^{11}_0,1)-\nu(\pi^{11}_0,7))$, 
$\mu(\pi^{11}_0,1)>\tfrac{4}{5}(\nu(\pi^{11}_0,1)-\nu(\pi^{11}_0,7))$,
as desired. 

Once again, it is straightforward to verify that 
we have a $1$--$1$ map from the partitions generated by
$\frac{1}{(q^3,q^5;q^8)_L}$
to the partitions generated by 
$\frac{1}{(q,q^7;q^8)_L}$,
which is not in general onto.  
This gives Theorem \ref{thm3} below. 

We also note the inequality \eqn{71}
is strict for $L\ge1$ and $n\ge 7$. The proof is similar to the $L=1$
case. 
We observe in each case that
\beq
\nu(\plz,1)\ge\nu(\plz,7).
\mylabel{eq:73b}
\eeq
As before, the partition $(7^m,1^r)$ does not satisfy \eqn{73b} when $m>r$,
and the remaining cases can be dealt with as before. 
And so we have proved the following

\begin{theorem}
\label{thm3}   
Let $\Atwid_{L,0}(n)$ denote the number of partitions $\tplz$ of $n$ into parts
$\equiv\pm 1\pmod{8}$, such that
the largest part $\le 8L-1$, $\nu(\tplz,1)\ge\nu(\tplz,7)$ and such that either
\begin{enumerate}
\item[(i)]
$3\mid (\nu(\tplz,1)-\nu(\tplz,7))$, and
$\mu(\tplz,1) \ge \tfrac{2}{3}(\nu(\tplz,1)-\nu(\tplz,7))$, or
\item[(ii)]
 $5\mid (\nu(\tplz,1)-\nu(\tplz,7))$, 
$3\nmid (\nu(\tplz,1)-\nu(\tplz,7))$, and 
$\mu(\tplz,1) \ge \tfrac{4}{5}(\nu(\tplz,1)-\nu(\tplz,7))$, or
\item[(iii)]
$\nu(\tplz,1)-\nu(\tplz,7)\equiv -8\pmod{15}$, $\mu(\tplz,7)>0$ and
$\mu(\tplz,1) \ge \tfrac{1}{5}(4(\nu(\tplz,1)-\nu(\tplz,7))-3)$.
\end{enumerate}
Then
\beq
\Atwid_{L,0}(n) = A_{L,1}(n),
\mylabel{eq:74}
\eeq
where $A_{L,1}(n)$ is defined in Theorem \ref{thm2}.
\end{theorem}

Theorems \ref{thm2} and \ref{thm3} with $n=19$ and $L=3$ are illustrated
below in Table \ref{table1}. The four partitions counted by $A_{3,1}(19)$
are listed in the first column. The first four partitions in the second
column are the corresponding images of our map, and are the partitions counted
by $\Atwid_{3,0}(19)$. The $8$ partitions counted $A_{3,0}(19)$ are
listed in the second column.

\begin{table}[ht]
\beqs
\renewcommand\arraystretch{1.1}
\begin{array}{lll}
(19^1) &\longrightarrow &(1^2,17^1) \\
(3^1,5^1,11^1) &\longrightarrow &(1^3,7^1,9^1) \\
(3^2,13^1) &\longrightarrow &(1^4,15^1) \\
(3^3,5^2) &\longrightarrow &(1^5,7^2) \\
          &                &(1^{19}) \\
          &                &(1^{12},7^1) \\
          &                &(1^{10},9) \\
          &                &(1^{1},9^2)\
\end{array}
\eeqs
\caption{Partitions of $n=19$, enumerated by $A_{3,1}$, $\Atwid_{3,0}$ and $A_{3,0}$}
\label{table1}
\end{table}

We would like to emphasise that the technique developed in this
sections is by no means limited to
\beqs
\frac{1}{(q,q^7;q^8)_L}
-
\frac{1}{(q^3,q^5;q^8)_L}.
\eeqs
In a very similar fashion we can prove the following theorems:
\begin{theorem}
\label{thm4}
Suppose $L>0$, and $1<r<m-1$. Then the
coefficients in the $q$-expansion of the difference of the two finite products
\beqs
\frac{1}{(q,q^{m-1};q^m)_L}-
\frac{1}{(q^r,q^{m-r};q^m)_L} 
\eeqs
are all nonnegative,
if and only if $r \nmid (m-r)$ and $(m-r)\nmid r$. 
\end{theorem}

\begin{theorem}
\label{thm5}   
Let $A^L_{m,r}(n)$ denote the number of partitions of $n$ into parts
$\equiv\pm r\pmod{m}$, with largest part $\le \max(Lm-r,Lm+r-m)$.
Let $\tilde A^L_{m,1}(n)$ denote the number of partitions $\tplm$ of $n$, counted 
by $A^L_{m,1}(n)$ subject to additional conditions that $\nu(\tplm,1)\ge\nu(\tplm,m-1)$, and
\begin{enumerate}
\item[(i)]
$r\mid(\nu(\tplm,1)-\nu(\tplm,m-1))$, and
$\mu(\tplm,1) \ge \tfrac{r-1}{r}(\nu(\tplm,1)-\nu(\tplm,m-1))$, or
\item[(ii)]
$(m-r)\mid (\nu(\tplm,1)-\nu(\tplm,m-1))$, 
$r\nmid(\nu(\tplm,1)-\nu(\tplm,m-1))$, and 
$\mu(\tplm,1) \ge \tfrac{m-r-1}{m-r}(\nu(\tplm,1)-\nu(\tplm,m-1))$, or
\item[(iii)]
$\nu(\tplm,1)-\nu(\tplm,m-1)\equiv -m\pmod{\frac{r(m-r)}{\gcd(r,m-r)}}$, $\mu(\tplm,m-1)>0$ and
$\mu(\tplm,1) \ge \tfrac{1}{m-r}((m-r-1)(\nu(\tplm,1)-\nu(\tplm,m-1))-r)$.
\end{enumerate}
Here, $\nu(\tplm,i)$ and $\mu(\tplm,i)$ denote the number of parts of $\tplm$ congruent to $i\pmod m$ 
and equal to $i$, respectively.\\
Then 
$$
\tilde A^L_{m,1}(n)=A^L_{m,r}(n),
$$ 
provided $0<r<m$, $r\nmid (m-r)$ and $(m-r)\nmid r$.
\end{theorem}

We plan to study more general partition
inequalities in a later paper.

Finally, we offer a prize of 500 $\yen$ for an elementary proof of the 
following conjectures:
\begin{align*}
c_{2n} &<  
\frac{\sqrt{\frac{13}{6}}}{32\sin(\pi/8)} \frac{1}{n}
\exp\{\tfrac{\pi}{2} \sqrt{\tfrac{13n}{3}}\},
\quad \mbox{for $n\ge1$, and}\\
c_{2n+1} &> 
\frac{\sqrt{\frac{13}{6}}}{32\cos(\pi/8)} \frac{1}{n}
\exp\{\tfrac{\pi}{2} \sqrt{\tfrac{13n}{3}}\},
\quad \mbox{for $n\ge2$.}
\end{align*}

\vskip 10pt
\noindent
{\it Acknowledgements}.
We are grateful to Holly Swisher for sending us her preprint \cite{Sw}
prior to publication. We would also like to thank George Andrews,
Krishna Alladi and Hamza Yesilyurt for helpful discussions.

\vskip 10pt
\noindent
{\it Note Added}.
After this manuscript was submitted for publication, Dennis Eichhorn 
brought to our attention the paper by Kevin Kadell, JCT A86, no.2 (1999), pp. 390-394.
In this paper the special case of our Theorem \ref{thm5} with $m=5$, $r=2$ and $L\rightarrow\infty$
is proven.

\bibliographystyle{amsplain}


\end{document}